\documentclass[12pt,a4paper]{article}
\usepackage{amsmath,amssymb,amsthm,amsfonts,latexsym}
\usepackage{texdraw}
\usepackage{graphicx}

\newtheorem{theorem}{Theorem}
\newtheorem{lemma}[theorem]{Lemma}
\newtheorem{corollary}[theorem]{Corollary}
\newtheorem{proposition}[theorem]{Proposition}

\theoremstyle{definition}

\theoremstyle{remark}

\newcommand{\z}{\mathbb{Z}}

\newcommand{\fd}{\mathbb{F}}
\newcommand{\dis}{\displaystyle}

\begin{document}
\title{Ternary Codes Associated with Symplectic Groups and Power Moments of Kloosterman Sums with Square Arguments}
\author{{DAE SAN KIM \footnote{Email address : dskim@sogang.ac.kr} and JI HYUN KIM \footnote{Email address : mip97@hanmail.net}} \\ \small{{Department of Mathematics, Sogang University, Seoul 121-742, South Korea}}}
\date{}
\maketitle

Abstract -- In this paper, we construct two ternary linear codes
associated with the symplectic groups $Sp(2,q)$ and $Sp(4,q)$.
Here $q$ is a power of three. Then we obtain recursive formulas
for the power moments of Kloosterman sums with square arguments
and for the even power moments of those in terms of the
frequencies of weights in the codes. This is done via Pless power
moment identity and by utilizing
the explicit expressions of ``Gauss sums" for the symplectic  groups $Sp(2n,q)$. \\

Key words -- ternary linear code, power moment, Kloosterman sum,
square argument, Pless power moment identity, Gauss sum,
symplectic group, weight distribution.\\

MSC 2000: 11T23, 20G40, 94B05.

\footnote{This work was supported by National Research Foundation of
Korea Grant funded by the Korean Government 2009-0072514.}

\section{Introduction}

 Let  $\psi$ be a nontrivial additive character of the finite field $\fd_q$
with $q=p^r$ elements ($p$ a prime). Then the Kloosterman sum
$K(\psi;a)$(\cite{R}) is defined by
\[
K(\psi;a) = \sum_{a \in \fd_q^*} \psi(\alpha + a\alpha^{-1}) (a
\in \fd_q^*).
\]
The Kloosterman sum was introduced in 1926(\cite{HD}) to give an
estimate for the Fourier coefficients of modular forms.

For each nonnegative integer $h$,  we denote by $MK(\psi)^h$ the
$h$-th moment of the Kloosterman sum $K(\psi;a)$, i.e.,
\[
MK(\psi)^h = \sum_{a \in \fd_q^*} K(\psi;a)^h.
\]
If $\psi = \lambda$ is the canonical additive character of
$\fd_q$, then $MK(\lambda)^h$ will be simply denoted by $MK^h$.

Also, we introduce an incomplete power moments of Kloosterman
sums. Namely, for every nonnegative integer $h$, and $\psi$ as
before, we define
\begin{equation}
SK(\psi)^h = \sum_{a \in \fd_q^*,~a~square} K(\psi;a)^h,
\end{equation}
which is called the $h$-th moment of Kloosterman sums with
``square arguments.'' If  $\psi = \lambda$ is the canonical
additive character of $\fd_q$, then $SK(\lambda)^h$ will be
denoted by $SK^h$, for brevity.

Explicit computations on power moments of Kloosterman sums were
initiated in the paper \cite{H} of Sali$\acute{e}$ in 1931, where it
is shown that for any odd prime $q$,
\[
MK^h = q^2 M_{h-1} - (q-1)^{h-1} + 2(-1)^{h-1} (h \geq 1).
\]
Here $M_0 = 0$, and for $h \in \z_{>0}$,
\[
M_h = \mid \{ (\alpha_1,\ldots,\alpha_h) \in (\fd_q^*)^h \mid
\sum_{j=1}^h \alpha_j = 1 = \sum_{j=1}^h \alpha_j^{-1} \} \mid.
\]
For  $q=p$ odd prime, Sali$\acute{e}$ obtained  $MK^1, MK^2, MK^3,
MK^4$ in \cite{H} by determining $M_1, M_2, M_3$. $MK^5$ can be
expressed in terms of the $p$-th eigenvalue for a weight 3 newform
on $\Gamma_0(15)$(cf. \cite{RL}, \cite{C}). $MK^6$ can be
expressed in terms of the $p$-th eigenvalue for a weight 4 newform
on $\Gamma_0(6)$(cf. \cite{K}). Also, based on numerical evidence,
in \cite{RJ} Evans was led to propose a conjecture which expresses
$MK^7$ in terms of Hecke eigenvalues for a weight 3 newform on
$\Gamma_0(525)$ with quartic nebentypus of conductor 105.

Assume now that $q=3^r$. Recently, Moisio was able to find explicit
expressions of $MK^h$, for $h \leq 10$(cf.\cite{M}). This was done,
via Pless power moment identity, by connecting moments of
Kloosterman sums and the frequencies of weights in the ternary Melas
code of length $q-1$, which were known by the work of  Geer, Schoof
and Vlugt in \cite{GV}. In \cite{D3},  two infinite families of
ternary linear codes associated with double cosets in the symplectic
group $Sp(2n,q)$ were constructed in order to generate infinite
families of recursive formulas for the power moments of Kloosterman
sums with square arguments and for the even power moments of those
in terms of the frequencies of weights in those codes.

In this paper, we will be able to produce two recursive formulas
generating power moments of Kloosterman sums with square arguments
over finite fields of characteristic three. To do that, we will
construct two ternary linear codes $C(Sp(2,q))$ and $C(Sp(4,q))$,
respectively associated with the symplectic groups $Sp(2,q)$ and
$Sp(4,q)$, and express those power moments in terms of the
frequencies of weights in each code. Then, thanks to our previous
results on the explicit expressions of ``Gauss sums" for the
symplectic  groups $Sp(2n,q)$ \cite{Kim1}, we can express the weight
of each codeword in the duals of the codes in terms of Kloosterman
sums with square arguments. Then our formulas will follow
immediately from the Pless power moment identity (cf. (29)).

Theorem 1 in the following(cf. (4), (5), (7), (8)) is the main
result of this paper. Henceforth, we agree that, for nonnegative
integers $a,b,c,$
\begin{equation}
{c \choose a,b} = \frac{c!}{a!b!(c-a-b)!}, ~if ~~a+b \leq c,
\end{equation}
and
\begin{equation}
{c \choose a,b} = 0, ~if ~~a+b > c.
\end{equation}

\begin{theorem}
Let $q=3^r$. Then we have the following.\\
(a) For $h=1,2,\ldots,$
\begin{multline}
SK^h = \sum_{j=0}^{h-1} (-1)^{h+j+1} {h \choose j}
(q^2-1)^{h-j}SK^j \\
+ q^{1-h} \sum_{j=0}^{min \{N_1,h\}} (-1)^{h+j} C_{1,j}
\sum_{t=j}^h t! S(h,t) 3^{h-t} 2^{t-h-j-1} {N_1-j \choose N_1-t},
\end{multline}
where $N_1 = |Sp(2,q)| = q(q^2-1)$, and $\{C_{1,j}\}_{j=0}^{N_1}$
is the weight distribution of the ternary linear code $C(Sp(2,q))$
given by

\begin{multline}
C_{1,j} = \sum {q^2 \choose \nu_1,\mu_1} {q^2 \choose \nu_{-1},
\mu_{-1}}\prod_{\beta^2-1 \neq 0 ~square} {q^2+q \choose
\nu_\beta, \mu_\beta}\\
\times \prod_{\beta^2-1 ~nonsquare} {q^2-q \choose \nu_\beta,
\mu_\beta} ~(j=0,\cdots,N_1).~~~~~~~~~~~~~~~~~~~~~~~~~~~
\end{multline}

Here the sum is over all the sets of nonnegative integers
$\{\nu_\beta\}_{\beta \in \fd_q}$ and $\{\mu_\beta\}_{\beta \in
\fd_q}$ satisfying $\dis\sum_{\beta \in \fd_q} \nu_\beta +
\dis\sum_{\beta \in \fd_q} \mu_\beta = j$ and $\dis\sum_{\beta \in
\fd_q} \nu_\beta \beta = \dis\sum_{\beta \in \fd_q} \mu_\beta
\beta$. In addition, $S(h,t)$ is the Stirling number of the second
kind defined by
\begin{equation}
S(h,t) = \frac{1}{t!}\sum_{j=0}^t (-1)^{t-j} {t \choose j} j^h.
\end{equation}\\
(b) For $h=1,2,\ldots,$
\begin{multline}
SK^{2h} = \sum_{j=0}^{h-1} (-1)^{h+j+1} {h \choose j}
(q^6-q^4-q^3-q^2+q+1)^{h-j} SK^{2j}\\
+ q^{1-4h} \sum_{j=0}^{min \{N_2,h\}} (-1)^{h+j} C_{2,j}
\sum_{t=j}^h t! S(h,t) 3^{h-t} 2^{t-h-j-1} {N_2-j \choose N_2-t},
\end{multline}
where $N_2 = |Sp(4,q)| = q^4(q^2-1)(q^4-1)$, and
$\{C_{2,j}\}_{j=0}^{N_2}$ is the weight distribution of the
ternary linear code $C(Sp(4,q))$ given by
\begin{multline}
C_{2,j} = \sum {q^4(\delta(2,q;0)+q^5-q^2-3q+3) \choose \nu_0,
\mu_0} \\\times \prod_{\beta \in \fd_q^*}
{q^4(\delta(2,q;\beta)+q^5-q^3-q^2-2q+3) \choose \nu_\beta,
\mu_\beta}(j=0,\ldots,N_2).
\end{multline}
Here the sum is over all the sets of nonnegative integers
$\{\nu_\beta\}_{\beta \in \fd_q}$ and $\{\mu_\beta\}_{\beta \in
\fd_q}$ satisfying $\dis\sum_{\beta \in \fd_q} \nu_\beta +
\dis\sum_{\beta \in \fd_q} \mu_\beta = j$ and $\dis\sum_{\beta \in
\fd_q} \nu_\beta \beta = \dis\sum_{\beta \in \fd_q} \mu_\beta
\beta$, and, for every $\beta \in \fd_q$, $\delta(2,q;\beta) =
\mid \{(\alpha_1,\alpha_2) \in (\fd_q^*)^2 \mid
\alpha_1+\alpha_1^{-1}+\alpha_2+\alpha_2^{-1} = \beta \}\mid$.
\end{theorem}

\section{$Sp(2n,q)$}

For more details about this section, one is referred to the paper
\cite{Kim1}. Throughout this paper, the following notations will
be used:
\begin{itemize}
 \item [] $q = 3^r$ ($r \in \z_{>0}$),
 \item [] $\fd_q$ = the finite field with $q$ elements,
 \item [] $Tr A$ = the trace of $A$ for a square matrix $A$,
 \item [] ${}^tB$ = the transpose of $B$ for any matrix $B$.
\end{itemize}\
The symplectic group $Sp(2n,q)$ over the field $\fd_q$ is defined
as:
\[
Sp(2n,q) = \{ w \in GL(2n,q) \mid {}^twJw = J \},
\]
with
\[
J = \left[%
\begin{array}{cc}
  0 & 1_n \\
  -1_n & 0 \\
\end{array}%
\right].
\]
Let $P=P(2n,q)$ be the maximal parabolic subgroup of $Sp(2n,q)$
defined by:
\[
P(2n,q) = \left\{\left[%
\begin{array}{cc}
  A & 0 \\
  0 & {}^tA^{-1} \\
\end{array}%
\right] \left[
\begin{array}{cc}
  1_n & B \\
  0 & 1_n \\
\end{array}%
\right] \mid A \in GL(n,q), {}^tB = B \right\}.
\]
Then, with respect to $P=P(2n,q)$, the Bruhat decomposition of
$Sp(2n,q)$ is given by
\begin{equation}
Sp(2n,q) = \coprod_{r=0}^n P \sigma_r P,
\end{equation}
where
\[
\sigma_r = \left[%
\begin{array}{cccc}
  0 & 0 & 1_r & 0 \\
  0 & 1_{n-r} & 0 & 0 \\
  -1_r & 0 & 0 & 0 \\
  0 & 0 & 0 & 1_{n-r} \\
\end{array}%
\right] \in Sp(2n,q).
\]
Put, for each $r$ with $0 \leq r \leq n$,
\[
A_r = \{w \in P(2n,q) \mid \sigma_r w \sigma_r^{-1} \in P(2n,q)
\}.
\]
Expressing $Sp(2n,q)$ as a disjoint union of right cosets of
$P=P(2n,q)$, the Bruhat decomposition in (9) can be written as
\begin{equation}
Sp(2n,q) = \coprod_{r=0}^n P\sigma_r(A_r\setminus P).
\end{equation}
The order of the general linear group $GL(n,q)$ is given by
\[
g_n = \prod_{j=0}^{n-1} (q^n-q^j) = q^{{n \choose 2}}
\prod_{j=1}^n (q^j-1).
\]
For integers  $n,r$ with $0 \leq r \leq n$, the $q$-binomial
coefficients are defined as:
\[
\left[ \substack{n \\ r}
 \right]_q = \prod_{j=0}^{r-1} (q^{n-j} - 1)/(q^{r-j} - 1).
\]
Then, for integers  $n,r$ with $0 \leq r \leq n$, we have
\begin{equation}
\frac{g_n}{g_{n-r}g_r} = q^{r(n-r)}\left[ \substack{n \\ r}
 \right]_q.
\end{equation}
In \cite{Kim1}, it is shown that
\begin{equation}
|A_r| = g_r g_{n-r} q^{{n+1 \choose 2}} q^{r(2n-3r-1)/2}.
\end{equation}
Also, it is immediate to see that
\begin{equation}
|P(2n,q)| = q^{{n+1 \choose 2}} g_n.
\end{equation}
So, from (11)-(13), we get
\begin{equation}
|A_r \setminus P(2n,q)| = q^{{r+1 \choose 2}}\left[ \substack{n \\
r} \right]_q,
\end{equation}
and
\begin{equation}
|P(2n,q)|^2 |A_r|^{-1} = q^{n^2} \left[ \substack{n \\ r}
 \right]_q q^{{r \choose 2}} q^r \prod_{j=1}^n (q^j-1).
\end{equation}
Also, from (10) and (15), we have
\begin{equation}
|Sp(2n,q)| = \sum_{r=0}^n |P(2n,q)|^2|A_r|^{-1} = q^{n^2}
\prod_{j=1}^n (q^{2j}-1),
\end{equation}
where one can apply the following $q$-binomial theorem with
$x=-q$:
\[
\sum_{r=0}^n \left[ \substack{n \\ r}
 \right]_q (-1)^r q^{{r \choose 2}} x^r = (x;q)_n,
\]
with $(x;q)_n = (1-x)(1-qx)\cdots(1-q^{n-1}x)$.

\section{Gauss sums for $Sp(2n,q)$}

The following notations will be employed throughout this paper.
\begin{itemize}
\item [] $tr(x) = x + x^3 + \cdots + x^{3^{r-1}}$ the trace
    function $\fd_q \rightarrow \fd_3$,
    \item [] $\lambda_0(x) = e^{2\pi i x/3}$ the canonical additive character of
    $\fd_3$,
    \item [] $\lambda(x) = e^{2\pi i tr(x)/3}$ the canonical additive character of
    $\fd_q$.
\end{itemize}
Then any nontrivial additive character $\psi$ of $\fd_q$ is given
by $\psi(x) = \lambda(ax)$, for a unique $a \in \fd_q^*$.

For any nontrivial additive character $\psi$ of $\fd_q$ and $a \in
\fd_q^*$, the Kloosterman sum  $K_{GL(t,q)}(\psi;a)$ for $GL(t,q)$
is defined as
\[
K_{GL(t,q)}(\psi;a) = \sum_{w \in GL(t,q)} \psi(Trw + aTrw^{-1}).
\]
Notice that, for $t=1$, $K_{GL(1,q)}(\psi;a)$ denotes the
Kloosterman sum $K(\psi;a)$.

In \cite{Kim1}, it is shown that $K_{GL(t,q)}(\psi;a)$ satisfies
the following recursive relation: for integers $t \geq 2$, $a \in
\fd_q^*$,
\begin{equation}
K_{GL(t,q)}(\psi;a) = q^{t-1} K_{GL(t-1,q)}(\psi;a)K(\psi;a) +
q^{2t-2}(q^{t-1}-1) K_{GL(t-2,q)}(\psi;a),
\end{equation}
where we understand that $K_{GL(0,q)}(\psi;a)=1$. From (17), in
\cite{Kim1} an explicit expression of the Kloosterman sum for
$GL(t,q)$ was derived.

\begin{theorem}[\cite{Kim1}] For integers $t \geq 1$, and $a \in \fd_q^*$,
the Kloosterman sum $K_{GL(t,q)}(\psi;a)$ is given by
\[
K_{GL(t,q)}(\psi;a) = q^{(t-2)(t+1)/2} \sum_{l=1}^{[(t+2)/2]}
q^lK(\psi;a)^{t+2-2l} \sum \prod_{\nu=1}^{l-1} (q^{j_\nu-2\nu}-1),
\]
where $K(\psi;a)$ is the Kloosterman sum and the inner sum is over
all integers $j_1,\ldots,j_{l-1}$ satisfying $2l-1 \leq j_{l-1}
\leq j_{l-2} \leq \cdots\leq j_1 \leq t+1$. Here we agree that the
inner sum is $1$ for $l=1$.
\end{theorem}
In Section 5 of \cite{Kim1}, it is shown that the Gauss sum for
$Sp(2n,q)$ is given by:
\begin{align*}
\sum_{w \in Sp(2n,q)} \psi(Trw) &= \sum_{r=0}^n |A_r\setminus
P|\sum_{w \in P} \psi(Trw\sigma_r)\\
&= q^{{n+1 \choose 2}} \sum_{r=0}^n |A_r \setminus P| q^{r(n-r)}
a_r K_{GL(n-r,q)}(\psi;1).
\end{align*}
Here  $\psi$ is any nontrivial additive character of $\fd_q$,
$a_0=1$, and, for $r \in \z_{>0}$, $a_r$ denotes the number of all
$r \times r$ nonsingular alternating matrices over $\fd_q$, which
is given by
\[
a_r = \left\{%
\begin{array}{ll}
    0, & \hbox{if $r$ is odd,} \\
    q^{\frac{r}{2}(\frac{r}{2}-1)}\prod_{j=1}^{\frac{r}{2}}(q^{2j-1}-1), & \hbox{if $r$ is even.} \\
\end{array}%
\right.
\]
(cf. \cite{Kim1},  Proposition 5.1). So
\[
\sum_{w \in Sp(2n,q)} \psi(Trw) = q^{{n+1 \choose 2}} \sum_{0 \leq
r \leq n,~r~even} q^{rn-\frac{1}{4}r^2}\left[ \substack{n \\
r} \right]_q \prod_{j=1}^{r/2}(q^{2j-1}-1) K_{GL(n-r,q)}(\psi;1).
\]

For our purposes, we only need the following two expressions of
the Gauss sums for $Sp(2,q)$ and $Sp(4,q)$. So we state them
separately as a theorem. Also, for the ease of notations, we
introduce
\[
G_1(q) = Sp(2,q),~~~ G_2(q) = Sp(4,q).
\]

\begin{theorem} Let $\psi$ be any nontrivial additive character of $\fd_q$. Then we have
\begin{align*}
&\sum_{w \in G_1(q)} \psi(Trw) = qK(\psi;1),\\
&\sum_{w \in G_2(q)} \psi(Trw) = q^4(K(\psi;1)^2+q^3-q).
\end{align*}
\end{theorem}
The next corollary follows from Theorems 4  and simple change of
variables.

\begin{corollary}Let  $\lambda$ be the canonical additive character of $\fd_q$,
and let $a \in \fd_q^*$. Then we have
\begin{align}
& \sum_{w \in G_1(q)} \lambda(aTrw) = qK(\lambda;a^2),\\
& \sum_{w \in G_2(q)} \lambda(aTrw) = q^4(K(\lambda;a^2)^2+q^3-q).
\end{align}
\end{corollary}

\begin{proposition}\emph{\textbf{([5, (5.3--5)])}} Let  $\lambda$ be the canonical additive character of $\fd_q$, $m \in \z_{\geq 0}$, $\beta \in \fd_q$.
Then
\begin{equation}
\sum_{a \in \fd_q^*} \lambda(-a\beta) K(\lambda;a^2)^m = q
\delta(m,q;\beta)-(q-1)^m,
\end{equation}
where, for $m \geq 1$,
\begin{equation}
\delta(m,q;\beta) = \mid\{(\alpha_1,\ldots,\alpha_m) \in
(\fd_q^*)^m \mid
\alpha_1+\alpha_1^{-1}+\cdots+\alpha_m+\alpha_m^{-1} = \beta
\}\mid,
\end{equation}
and
\[
\delta(0,q;\beta) = \left\{%
\begin{array}{ll}
    1, & \hbox{if $\beta=0$,} \\
    0, & \hbox{otherwise.} \\
\end{array}%
\right.
\]
\end{proposition}

Let $G(q)$ be one of finite classical groups over $\fd_q$. Then we
put, for each $\beta \in \fd_q$,
\[
N_{G(q)}(\beta) = \mid\{w \in G(q) \mid Tr(w) = \beta \}\mid.
\]
Then it is easy to see that
\begin{equation}
qN_{G(q)}(\beta) = |G(q)| + \sum_{a \in \fd_q^*} \lambda(-a\beta)
\sum_{w \in G(q)} \lambda(aTrw).
\end{equation}
For brevity, we write
\begin{equation}
n_1(\beta) = N_{G_1(q)}(\beta), n_2(\beta) = N_{G_2(q)}(\beta).
\end{equation}
Using (16), (18)-(20), and (22), one derives the following
proposition. Here one notes that
\begin{align*}
\delta(1,q;\beta) &= \mid\{x \in \fd_q \mid x^2-\beta x +1 = 0
\}\mid\\
&= \left\{%
\begin{array}{lll}
    2, & \hbox{if $\beta^2-1 \neq 0$ is a square,} \\
    1, & \hbox{if $\beta = \pm 1$,} \\
    0, & \hbox{if $\beta^2-1$ is a nonsquare.} \\
\end{array}%
\right.
\end{align*}

\begin{proposition} With  $n_1(\beta), n_2(\beta)$ as in (23) and $\delta(m,q;\beta)$ as in (21),  we have:
\begin{equation}
n_1(\beta) = q \delta(1,q;\beta)+ q^2-q = \left\{%
\begin{array}{ll}
    q^2+q, & \hbox{if $\beta^2-1 \neq 0$ is a square,} \\
    q^2, & \hbox{if $\beta = \pm 1$,} \\
    q^2-q, & \hbox{if $\beta^2-1$ is a nonsquare.} \\
\end{array}%
\right.
\end{equation}
\begin{equation}
n_2(\beta) = \left\{%
\begin{array}{ll}
    q^4(\delta(2,q;0)+q^5-q^2-3q+3), & \hbox{if $\beta = 0$,} \\
    q^4(\delta(2,q;\beta) + q^5-q^3-q^2-2q+3), & \hbox{if $\beta \neq 0$.} \\
\end{array}%
\right.~~~~~~~~~~~~
\end{equation}
\end{proposition}

\begin{proposition}
$Tr : Sp(2n,q) \rightarrow \fd_q$ is surjective.
\end{proposition}
\begin{proof}
Under the map $Tr$, for any $\alpha \in \fd_q$,
\[
Sp(2n,q) \ni \left[%
\begin{array}{cccc}
  \alpha & 0 & 1 & 0 \\
  0 & 0 & 0 & 1_{n-1} \\
  -1 & 0 & 0 & 0 \\
  0 & -1_{n-1} & 0 & 0 \\
\end{array}%
\right] \mapsto \alpha.
\]
\end{proof}

\section{Construction of codes}

Let
\begin{equation}
N_1 = |G_1(q)| = q(q^2-1), ~~N_2 = |G_2(q)| = q^4(q^2-1)(q^4-1).
\end{equation}
Here we will construct two ternary linear codes  $C(G_1(q))$ of
length $N_1$ and $C(G_2(q))$ of length $N_2$, respectively
associated with the symplectic groups $G_1(q)$ and $G_2(q)$.

By abuse of notations, for $i=1,2,$ let $g_1,g_2,\ldots,g_{N_i}$
be a fixed ordering of the elements in the group $G_i(q)$. Also,
for $i=1,2,$  we put
\[
v_i = (Trg_1,Trg_2,\ldots,Trg_{N_i}) \in \fd_q^{N_i}.
\]
Then, for $i=1,2,$ the ternary linear code  $C(G_i(q))$ is defined
as
\begin{equation}
C(G_i(q)) = \{ u \in \fd_3^{N_i} \mid u\cdot v_i = 0 \},
\end{equation}
where the dot denotes the usual inner product in $\fd_q^{N_i}$.

The following Delsarte's theorem is well-known.
\begin{theorem}[\cite{FJ}]  Let  $B$ be a linear code over $\fd_q$.  Then
\[
(B|_{\fd_3})^\bot = tr(B^\bot).
\]
\end{theorem}
In view of this theorem, the dual $C(G_i(q))^\bot$ ($i=1,2$) is
given by
\begin{equation}
C(G_i(q))^\bot = \{c_i(a) = (tr(aTrg_1),\ldots,tr(aTrg_{N_i}))
\mid a \in \fd_q \}.
\end{equation}

\begin{proposition}
 For every $q=3^r$, and $i=1,2,$ the map $\fd_q \rightarrow C(G_i(q))^\bot$ $(a \mapsto
 c_i(a))$
 is an $\fd_3$-linear isomorphism.
\end{proposition}
\begin{proof}
The map is clearly $\fd_3$-linear and surjective. Let $a$ be in
the kernel of the map. Then, in view of Proposition 7,
$tr(a\beta)=0$, for all $\beta \in \fd_q$. Since the trace
function $\fd_q \rightarrow \fd_3$  is surjective, $a=0$.
\end{proof}

\section{Recursive formulas for power moments of Kloosterman sums with square arguments}

In this section, we will be able to find, via Pless power moment
identity,  recursive formulas for the power moments of Kloosterman
sums with square arguments and even power moments of those with
square arguments  in terms of the frequencies of weights in
$C(G_i(q))$, for each $i=1,2.$

\begin{theorem}[Pless power moment identity]
Let $B$ be an $q$-ary $[n,k]$ code, and let $B_i$(resp.
$B_i^\bot$) denote the number of codewords of weight $i$ in
$B$(resp. in $B^\bot$). Then, for $h=0,1,2,\ldots,$
\begin{equation}
\sum_{j=0}^n j^h B_j = \sum_{j=0}^{min \{n,h\}} (-1)^j B_j^\bot
\sum_{t=j}^h t! S(h,t) q^{k-t} (q-1)^{t-j} {n-j \choose n-t},
\end{equation}
where $S(h,t)$ is the Stirling number of the second kind defined
in (6).
\end{theorem}

\begin{lemma}
Let $c_i(a) = (tr(aTrg_1),\ldots,tr(aTrg_{N_i})) \in
C(G_i(q))^\bot$, for $a \in \fd_q^*$, and $i=1,2.$ Then the
Hamming weight $w(c_i(a))$ can be expressed as follows:
\begin{align}
(a)~ w(c_1(a)) &= \frac{2}{3}q(q^2-1-K(\lambda;a^2)),~~~~~~~~~~~~~~~~~~~~~~~~~~~~~~~~~~~~~~\\
(b)~ w(c_2(a)) &= \frac{2}{3}q^4\{(q^2-1)(q^4-1)-(K(\lambda;a^2)^2
+ q^3-q) \}.
\end{align}
\end{lemma}
\begin{proof} For $i=1,2,$
\begin{align*}
w(c_i(a)) &= \sum_{j=1}^{N_i} (1-\frac{1}{3}\sum_{\alpha \in
\fd_3} \lambda_0(\alpha
tr(aTrg_j))~~~~~~~~~~~~~~~~~~~~~~\\
&= N_i - \frac{1}{3}\sum_{\alpha \in \fd_3} \sum_{w \in G_i(q)}
\lambda(\alpha a Trw)\\
&= \frac{2}{3}N_i - \frac{1}{3}\sum_{\alpha \in \fd_3^*} \sum_{w
\in G_i(q)} \lambda(\alpha a Trw).
\end{align*}
Our results now follow from (18), (19), and (26).
\end{proof}
Fix $i(i=1,2),$ and let $u=(u_1,\ldots,u_{N_i}) \in \fd_3^{N_i}$,
with $\nu_\beta$ 1's and $\mu_\beta$ 2's in the coordinate places
where $Tr(g_j)=\beta$, for each $\beta \in \fd_q$. Then we see
from the definition of the code $C(G_i(q))$(cf. (27)) that $u$ is
a codeword with weight $j$ if and only if $\dis\sum_{\beta \in
\fd_q} \nu_\beta + \dis\sum_{\beta \in \fd_q} \mu_\beta = j$ and
$\dis\sum_{\beta \in \fd_q} \nu_\beta \beta = \dis\sum_{\beta \in
\fd_q} \mu_\beta \beta$(an identity in $\fd_q$). Note that there
are $\dis \prod_{\beta \in \fd_q} {n_i(\beta) \choose \nu_\beta,
\mu_\beta}$ (cf. (2), (3)) many such codewords with weight $j$.
Now, we get the following formulas in (32) and (33), by using the
explicit values of  $n_i(\beta)$ in (24) and (25).

\begin{theorem}
Let $\{C_{i,j}\}_{j=0}^{N_i}$ be the weight distribution of
$C(G_i(q))$, for $i=1,2.$ Then
\begin{multline}
(a)~~ C_{1,j} = \sum {q^2 \choose \nu_1,\mu_1} {q^2 \choose
\nu_{-1}, \mu_{-1}}\prod_{\beta^2-1 \neq 0 ~square} {q^2+q \choose
\nu_\beta, \mu_\beta}\\
\times \prod_{\beta^2-1 ~nonsquare} {q^2-q \choose \nu_\beta,
\mu_\beta} ~(j=0,\ldots,N_1),~~~~~~~~~~~~~~~~~~~~
\end{multline}
where the sum is over all the sets of nonnegative integers
$\{\nu_\beta\}_{\beta \in \fd_q}$ and $\{\mu_\beta\}_{\beta \in
\fd_q}$ satisfying
\[\dis\sum_{\beta \in \fd_q} \nu_\beta +
\dis\sum_{\beta \in \fd_q} \mu_\beta = j ~~\text{and}~~
\dis\sum_{\beta \in \fd_q} \nu_\beta \beta = \dis\sum_{\beta \in
\fd_q} \mu_\beta \beta.
\]
\begin{multline}
(b)~~C_{2,j} = \sum {q^4(\delta(2,q;0)+q^5-q^2-3q+3) \choose
\nu_0, \mu_0} \\\times \prod_{\beta \in \fd_q^*}
{q^4(\delta(2,q;\beta)+q^5-q^3-q^2-2q+3) \choose \nu_\beta,
\mu_\beta}(j=0,\ldots,N_2),
\end{multline}
where the sum is over all the sets of nonnegative integers
$\{\nu_\beta\}_{\beta \in \fd_q}$ and $\{\mu_\beta\}_{\beta \in
\fd_q}$ satisfying
\[\dis\sum_{\beta \in \fd_q} \nu_\beta +
\dis\sum_{\beta \in \fd_q} \mu_\beta = j ~~\text{and}~~
\dis\sum_{\beta \in \fd_q} \nu_\beta \beta = \dis\sum_{\beta \in
\fd_q} \mu_\beta \beta,
\]
and, for every $\beta \in \fd_q$,
\[
\delta(2,q;\beta) = \mid\{(\alpha_1,\alpha_2) \in (\fd_q^*)^2 \mid
\alpha_1 + \alpha_1^{-1} + \alpha_2 + \alpha_2^{-1} = \beta
\}\mid.
\]
\end{theorem}

We now apply the Pless power moment identity in (29) to each
$C(G_i(q))^\bot$, for $i=1,2,$ in order to obtain the results in
Theorem 1(cf. (4), (7)) about recursive formulas.

Then the left hand side of the identity in (29) is equal to
\begin{equation}
\sum_{a \in \fd_q^*} w(c_i(a))^h,
\end{equation}
with the $w(c_i(a))$ in each case given by (30) and (31).\\

For $i=1,$  (34) is
\begin{equation}
\begin{split}
& (\frac{2q}{3})^h \sum_{a \in \fd_q^*} (q^2-1-K(\lambda;a^2))^h\\
&= (\frac{2q}{3})^h \sum_{a \in \fd_q^*}\sum_{j=0}^h (-1)^j {h
\choose j} (q^2-1)^{h-j} K(\lambda;a^2)^j\\
&= 2(\frac{2q}{3})^h \sum_{j=0}^h (-1)^j {h \choose j}
(q^2-1)^{h-j} SK^j.
\end{split}
\end{equation}
Similarly, for $i=2,$  (34) equals
\begin{equation}
2(\frac{2q^4}{3})^h \sum_{j=0}^h (-1)^j {h \choose j}
(q^6-q^4-q^3-q^2 + q +1)^{h-j} SK^{2j}.
\end{equation}
Here one has to separate the term corresponding to $j=h$  in (35)
and (36), and note $dim_{\fd_3}C(G_i(q))^\bot = r$.


\begin{thebibliography}{1}

\bibitem{RJ}
 R.J.Evans, ``Seventh power moments of Kloosterman sums,''
  \emph{Israel J. Math.}, to appear.

\bibitem{GV}
G. van der Geer, R. Schoof and M. van der Vlugt, ``Weight formulas
for ternary Melas codes,''
  \emph{Math. Comp.} 58(1992), 781--792.

\bibitem{K}
 K. Hulek, J. Spandaw, B. van Geemen, and D. van Straten, ``The modularity
of the Barth-Nieto quintic and its relatives,''\emph{Adv. Geom.} 1
(2001), 263-289.

\bibitem{Kim1}
D.S.Kim, ``Gauss sums for symplectic groups over a finite field,''
\emph{Mh. Math.} 126(1998), 55-71.

\bibitem{Kim2}
D.S.Kim, ``Exponential sums for symplectic groups and their
applications,'' \emph{Acta Arith.} 88(1999), 151-171.

\bibitem{D3}D. S. Kim,
``Infinite families of recursive formulas generating power moments
of ternary Kloosterman sums with square arguments arising from
symplectic groups,'' \emph{Adv. Math. Commun.} 3(2009), 167--178.

\bibitem{HD}
H.D. Kloosterman , ``On the representation of numbers in the form
$ax^2 + by^2 + cz^2 + dt^2$,''\emph{Acta Math.} 49(1926),
407--464.

\bibitem{R}
R. Lidl and H.Niederreiter, ``Finite fields,'' \emph{Encyclopedia
of Mathematics and Its Applications 20, Cambridge University
Pless, Cambridge}, 1987.

\bibitem{RL}
R. Livn$\acute{e}$, ``Motivic orthogonal two-dimensional
representations of $Gal(\overline{Q}/Q)$,'' \emph{Israel J. Math.}
92 (1995), 149-156.

\bibitem{FJ}
F.J. MacWilliams and N.J. A. Sloane, ``The theory of error
correcting codes,'' \emph{North-Holland, Amsterdam}, 1998.

\bibitem{M}
M. Moisio, ``On the moments of Kloosterman sums and fibre products
of Kloosterman curves,'' \emph{Finite Fields Appl.} 14(2008),
515-531.

\bibitem{C}
 C. Peters, J. Top, and M. van der Vlugt, ``The Hasse zeta function of a K3
surface related to the number of words of weight 5 in the Melas
codes,'' \emph{ J. Reine Angew. Math.} 432 (1992), 151-176.

\bibitem{H}
H. Sali\'{e}, ``Uber die Kloostermanschen Summen $s(u,v;q)$ ,''
\emph{Math. Z.}  34(1931), 91--109.

\end{thebibliography}
\end{document}